\documentclass[11pt]{amsart}
\usepackage[margin=1.0in]{geometry}
\usepackage{amsmath,amssymb,amsfonts,graphicx,color,fancyhdr,psfrag,comment,enumerate}
\usepackage[hyperpageref]{backref}
\usepackage[colorlinks=true, pdfstartview=FitV, linkcolor=blue, citecolor=blue, urlcolor=blue]{hyperref}

\usepackage[nodisplayskipstretch]{setspace}
\usepackage[utf8]{inputenc}

\usepackage{epstopdf}
\usepackage{mathrsfs}
\usepackage{epsfig}
\usepackage{hyperref}
\usepackage{ifthen}
\usepackage{float}
\usepackage{enumerate}
\usepackage{mathtools}
\usepackage[bottom]{footmisc}
\usepackage{tikz}
\usetikzlibrary{matrix,shapes,arrows,positioning,chains}
\usepackage{amssymb,amsmath,amsfonts}
\usepackage{bbm}
\usepackage[dvipsnames]{xcolor}
\usepackage[normalem]{ulem}

\numberwithin{equation}{section}

\DeclareMathOperator\supp{supp}
\newcommand{\R}{\mathbb R}
\newcommand{\N}{\mathbb N}

\newcommand{\Z}{\mathbb Z}

\newcommand{\A}{\mathcal A}

\newtheorem{theorem}{Theorem}[section]
\newtheorem{prop}[theorem]{Proposition}
\newtheorem{lemma}[theorem]{Lemma}

\newtheorem{remark}[theorem]{Remark}

\usepackage{tikz}
\usepackage{pgfplots}
\pgfplotsset{compat=1.13,
tick label style={color=white},
  label style={font=\small},
  legend style={font=\small}}

\usepackage{mathtools}

\allowdisplaybreaks

\newcommand{\NN}{\mathfrak N}
\newcommand{\MM}{\mathfrak M}

\begin{document}
	
	\title{Variation bounds for averages over nondegenerate curves}

	\begin{abstract}
	In this article, we prove sharp variational estimates for averaging operators over non-degenerate curves in $\R^3.$ The sufficiency follows from  dyadic localisation, Sobolev embedding, and local smoothing estimates, and the primary novelty of this article lies in obtaining sharp necessary conditions  for variational estimates.  
	\end{abstract}

\author[Ghosh]{Abhishek Ghosh}

\address[Abhishek Ghosh]{Department of Mathematics, Indian Institute of Technology Madras, Chennai, 600036, India.}

\email{abhi@iitm.ac.in, abhi170791@gmail.com}

\author[Shuin]{Kalachand Shuin}

\address[Kalachand Shuin]{Kerala School of Mathematics, Kozhikode-673571, Kerala, India.}
\email{kcshuin@ksom.res.in, shuin.k.c@gmail.com}

\maketitle
 \section{Introduction}

Variational estimates have gathered a lot of attention in the last decades due to its significance
in strengthening pointwise convergence results for sequence of averaging operators and it goes back to the fundamental work of L\'epingle in \cite{Lep} regarding martigle sequences. Variation-norm estimates for averaging operators over curved submanifolds is an important direction starting with the pioneering work of Jones, Seeger, and Wright \cite{JSW}. To clearly state their results, recall  that for a  family of complex numbers $\textbf{a}=\{a_t\}_{t>0},$ the $q$-variation of $\textbf{a}$ is defined as
$$
\|\textbf{a}\|_{v_q}=\sup_{L\in\mathbb{N}}\sup_{0<t_1<t_2<\cdots<t_{L}<\infty}
\left(\sum_{i=1}^{L-1}\left|a_{t_{i+1}}-a_{t_i}\right|^q\right)^{\frac1q},
$$
 for $q\in[1,\infty)$ and replacing the $l^{q}$-sum by a supremum when $q=\infty$. In \cite{JSW} the authors showed that variation norm operator associated with spherical averages, defined by 
\[
V_q(\{S_{t}f\})(x)=\big\|\{S_t f(x)\}_{t>0}\big\|_{v^q}, \quad\text{where,}\quad S_t f(x)=\int_{S^{n-1}}f(x-ty)\,d\sigma(y)
\]
and \(d\sigma\) denotes the normalized surface measure on \(S^{n-1}\), is bounded on $L^p(\mathbb R^n)$ for all $q>2$ if $\frac{n}{n-1}<p\leq 2n;$ and when $p>2n$ then $V_q(\{S_{t}\}): L^p(\R^n)\to L^p(\R^n)$ if $q>p/n.$ They further
showed that the $L^p$ boundedness fails at $q=2$, and also
whenever $q<\frac{p}{n}$ for $p>2n.$ Since $V_{r}(\{S_{t}f\})(x)\geq \sup_{t}|S_{t}f(x)-S_{t_{0}}f(x)|$ for a.e. $x\in \R^n,$ $t_{0}\in (0, \infty),$ these provides a significant strengthening of the spherical maximal theorem established by Stein \cite{SteinSph}, for $n>2$, and Bourgain \cite{Bourgain}, for $n=2$, stating that the spherical maximal operator $\mathcal{S}f=\sup_{t>0}|S_t f(x)|$
is bounded on $L^p(\mathbb R^n)$ if, and only if $p>\frac{n}{n-1}.$ More recently, Beltran, Oberlin, Roncal, et al. \cite{BOR} addressed the remaining endpoint $q=p/n$ in \cite{JSW}. In dimensions $n>2$, they proved that, for $p>2n,$
\[
V_{p/n}(\{S_{t}\})
:
L^{p,1}(\mathbb R^n)
\longrightarrow
L^{p,\infty}(\mathbb R^n).
\]
Importantly, they also proved sharp $L^p$ improving estimates for local variational operators defined by $V_q(\{S_{t}f\}; t\in [1, 2]):=\left\|\{S_t f(x)\}_{t\in [1, 2]}\right\|_{v_q}$ which strengthens $L^p$ improving estimates for the local spherical maximal function obtained in \cite{Schlag, SS}.  As a consequence, they further derive sparse operator bounds and weighted estimates in \cite{BOR}. In \cite{gensph} the authors obtained sharp variational estimates for the generalized spherical means, and very recently in \cite{finite} variational estimates for averages over finite type curves on the plane were obtained.

Motivated by these works, the primary goal of this article is to obtain sharp variational estimates for averages over non-degenerate curves on $\R^3.$ To provide context
and clearly state our main results, we start with defining our main objects.

Let $\gamma:I\to \R^n$ be a smooth non-degenerate curve, that is
\begin{eqnarray}\label{nondegenerate}
    |\det(\gamma'(s),\gamma''(s),\gamma'''(s), \cdots, \gamma^{n}(s))|\geq c_0,
\end{eqnarray}
for some positive constant $c_0$. For $t>0,$ we consider the following averaging operator and the associated maximal function
\begin{eqnarray}\label{average}
    A_tf(x):=\int_{I} f(x-t\gamma(s)) \chi(s)~ds,\quad Mf(x):=\sup_{t>0}|A_{t}f(x)|,
\end{eqnarray}
where $\chi$ is a bump function supported in the interior of $I$. In recent times, the study of averaging operators and maximal functions over non-degenerate curves has gained a lot of attraction and we start with the  seminal result proved independently by Beltran \textit{et al} in \cite{BeltranAJM} and by Ko--Lee-Oh in \cite{KoLeeOh} that the maximal operator $M$ is bounded from $L^p(\mathbb{R}^3)$ to itself if and only if $p>3.$ We also mention the work by Pramanik and Seeger in \cite{PramanikAJM} where they proved the estimate for large $p.$ A general strategy to prove maximal function estimate relies on proving local--smoothing estimates of the form 
\begin{align}\Big(\int^2_1 \int_{\R^3} |A_t P_{k} f(x)|^p ~dx dt\Big)^{1/p}\lesssim_{\epsilon} 2^{-k\sigma(p)} 2^{\epsilon k} \Vert f\Vert_{L^{p}(\R^3)},
\label{general:ls}\end{align}
and this goes back to the work in \cite{MSS} in the context of the circular averaging operators. Here $P_j$ denotes the smooth Littlewood--Paley projection localised to frequencies $|\xi|\sim 2^{k}.$ In \cite{BeltranAJM} the authors obtained \eqref{general:ls} for some $\sigma>\frac{1}{p}$ for $p>3$ consequently obtained the maximal estimate as an application of Sobolev embedding. An independent approach was taken by Ko--Lee--Oh in \cite{KoLeeOh}, and subsequently in \cite{Forummathpi} they achieved maximal function estimates over nondegenerate curves on $\R^n$ for $p>2n-4, n\geq 4,$ leaving the range $n<p\leq 2n-4$ open. A striking example in \cite{BH} showed that the resolution to the sharp local smoothing is not sufficient to prove sharp maximal function estimates when $n\geq 5,$ however, recently in \cite{Changkeon} Gan–Maldague–Oh first proved sharp local smoothing estimates for averages over nondegenerate curves in any dimension $n\geq 3,$ and subsequently in \cite{Changkeon1} Gan-Oh--Shin obtained the complete resolution to the maximal function estimate proving that $M:L^p(\R^n)\to L^p(\R^n)$ if and only if $p>n, n\geq 3.$ Motivated by these results, in this article, we obtain sharp variational estimates for averages over non-degenerate curves on $\R^3.$ Let $\mathcal{A}=\{A_t\}_{t>0}$  be the family of the averages defined by \eqref{average}, and define
$$
V_q(\mathcal{A}f)(x):=\left\|\{A_tf(x)\}_{t>0}\right\|_{v_q}.
$$ 
Also, $V_q(\mathcal{A}f; E):=\left\|\{A_tf(x)\}_{t\in E}\right\|_{v_q}.$ Now we state our main result. 

\begin{theorem}\label{mainresult}
    Let $\mathcal{A}=\{A_t\}_{t>0}$ defined in \eqref{average} with $\gamma$ satisfying \eqref{nondegenerate}. Let $3<p<\infty$. Then the following boundedness holds 
    \begin{eqnarray}\label{variationbound}
       \Vert V_q(\mathcal{A}f)\Vert_{L^{p}}\leq C(p, q) \Vert f\Vert_{L^p}
    \end{eqnarray}
    if $q>\max\{3,\frac{6p}{p+4},\frac{p}{2}\}$. Moreover, the above boundedness fails if $q<\max\{3,\frac{6p}{p+4},\frac{p}{2}\}$.
\end{theorem}
Since the variational seminorm $V_q$ increases as $q$ decreases, our results strengthen the maximal theorem obtained in \cite{BeltranAJM, KoLeeOh}. We conclude this section by a few words on the proofs and some preliminaries.

\subsection{Preliminaries}  For $j\in\mathbb Z$, define
\begin{align*}
V_{q,j}(\mathcal Af)(x)
=
\sup_{\substack{t_1<\cdots<t_N\\
[t_\ell,t_{\ell+1}]\subset[2^j,2^{j+1}]}}
\left(
\sum_{\ell=1}^{N-1}
\left|A_{t_{\ell+1}}f(x)-A_{t_\ell}f(x)\right|^q
\right)^{1/q}.
\label{defn}
\end{align*}
The variational operator is decomposed into the short $q$-variation and long variation defined by
\begin{align*}
S_{q}(\mathcal Af)(x)
=
\big(
\sum_{j\in\mathbb Z}
\big[V_{q,j}(\mathcal Af)(x)\big]^q
\big)^{1/q},\quad V_q^{\mathrm{long}}(\mathcal Af)(x)
:=
\sup_{N\in\mathbb N}
\sup_{\substack{t_1<\cdots<t_N\\ t_i\in\mathbb Z}}
\big(
\sum_{i=1}^{N-1}
\left|
A_{2^{t_{i+1}}}f(x)-A_{2^{t_i}}f(x)
\right|^q
\big)^{1/q},
\end{align*}
respectively. Consequently, for $q>2,$ we have $$V_{q}(\mathcal{A} f)(x)\lesssim S_{q}(\mathcal Af)(x)+ V_q^{\mathrm{long}}(\mathcal Af)(x).$$ We write $\widehat {A_{t}f}(\xi)=\widehat \mu_{t}(\xi) \widehat f(\xi),$ where $\widehat{\mu_{t}}(\xi):=\widehat{\mu}(t\xi)= \int e^{-i\langle t\gamma(s), \xi)\rangle}\,ds.$ By Van der Corput  lemma we get $$|\widehat{\mu_{t}}(\xi)|\lesssim (1+t|\xi|)^{-1/3},$$
which combined with the in \cite[Theorem~1.1]{JSW} concludes that the long variation satisfies the estimate
\begin{align*}
\|V_q^{\mathrm{long}}(\mathcal Af)\|_{L^p(\R^3)}\leq C \|f\|_{L^p(\R^3)}  
\end{align*}
for all $1<p<\infty$ and $q>2.$ Therefore, to prove the sufficiency in Theorem~\ref{mainresult} we will only focus on the short variation $S_{q}(\mathcal Af).$ To handle the short variation we require sharp local smoothing estimates obtained in \cite{Changkeon}. Define the Littlewood-Paley projections $\widehat{P_{k}f}({\xi})=\widehat{\beta_{k}(D)f}(\xi)=\beta_{k}(\xi) \widehat{f}(\xi),$ where $\beta_{k}(\xi):=\beta(|\xi|/2^{k})$ for $k\geq1$, and the smooth function $\beta$ is defined as $\beta(r):=\zeta(r)-\zeta(2r)$ and $\zeta$ is a non-negative, even, compactly supported smooth function with $\zeta=1$ on $[-1, 1]$ and $\zeta=0$ outside $[-2, 2]$. 
\begin{theorem}\cite[Theorem 1.1]{Changkeon}\label{LS:thm}
Let $\gamma$ be a smooth non-degenerate curve satisfying \eqref{nondegenerate}. Let $2\leq p<\infty$ and denote $\sigma(p):=\min{\{\frac{1}{3}, \frac{1}{3}(\frac{1}{2}+\frac{2}{p}), \frac{2}{p}\}}.$ Then for all $k\geq1,$ and for all $\epsilon>0$  we have
\[\Big(\int^2_1 \int_{\R^3} |A_tP_kf(x)|^p ~dx dt\Big)^{1/p}\lesssim_{\epsilon} 2^{-k\sigma(p)} 2^{\epsilon k} \Vert f\Vert_{L^{p}(\R^3)}.\]
\end{theorem}
We shall also use the Sobolev embedding: For $1\leq q<\infty,$
\begin{align}
\label{sobolev}
\|a\|_{v_q}
\leq C_{q} \|a\|_{L^q(\R^+)}^{1/q'}
\|a'\|_{L^q(\R^+)}^{1/q}.
\end{align}

A key novelty of this article lies in proving the necessary conditions in Theorem~\ref{mainresult}. Necessity of $q\geq 3,$ and $q\geq \frac{p}{2}$ in Theorem \ref{mainresult} follows from modifying examples from \cite{BOR} for the averages over non-degenerate curves. However, to show that the condition $q\geq \frac{6p}{p+4}$ is necessary, we construct examples motivated by \cite{BH, BeltranAdv} with crucial modifications pertaining to the computations related to the variational seminorm. In Section~\ref{sec1} we prove the necessity of $q\geq 3,$ and $q\geq \frac{p}{2}.$ Section~\ref{sec2} contains our main counterexample showing necessity of $q\geq \frac{6p}{p+4},$ and finally in Section~\ref{suff} we prove the sufficiency in Theorem~\ref{mainresult}.
 
 \section{Necessity of $q\geq 3,$ and $q\geq \frac{p}{2}$ in Theorem \ref{mainresult}}
 \label{sec1}
 We shall prove the necessary part of Theorem \ref{mainresult} in the following few propositions.
\begin{prop}
Let \eqref{variationbound} holds true for all $f\in L^p.$ Then we necessarily have  
 \begin{equation*}
\label{cond:1}
 q\geq \frac{p}{2}.   
\end{equation*}
\end{prop}
\begin{proof}
Throughout, \[ \gamma(s)=\left(s,\frac{s^2}{2},\frac{s^3}{6}\right),\qquad A_tf(x)=\int_{I} f(x-t\gamma(s))\chi(s)\,ds, \] where \(\chi\ge0\) is smooth and positive on a fixed interval around some \(s_0>0\). Define \[ t_j=\frac54+jL\delta,\qquad 1\le j\le N, \] where \(L\) is a sufficiently large fixed constant and \[ N\le c\delta^{-1}. \] Define \begin{align} \mathcal N_j(\delta) = \{y\in \R^3: |y-t_j\gamma(s)|\leq \delta\,\quad \text{for some}\,\, s\in I\}, \qquad f_\delta=\sum_{j=1}^N(-1)^j\mathbf1_{\mathcal N_j(\delta)}.\label{nbd:ex:1}\end{align} 

For $t\in [1, 2], s\in I,$ one can observe that 
\begin{equation}
\label{sep:ex:1}
 |t\gamma(s)-t'\gamma(s')| \gtrsim |t-t'|+|s-s'|.    
\end{equation}
Thus, choosing \(L\) large makes the neighbourhoods in \eqref{nbd:ex:1} pairwise disjoint. Moreover, $|\mathcal N_j(\delta)|\sim \delta^2.$ Consequently, \[ \|f_\delta\|_p\sim (N\delta^2)^{1/p}. \label{norm:ex:1} \] Now let \(|x|\le c_0\delta\), with \(c_0\) small. At time \(t_j\), \[ x-t_j\gamma(s)\in\mathcal N_j(\delta) \qquad(s\in I). \] Therefore \[ A_{t_j}f_\delta(x)=(-1)^j a_j(x), \qquad a_j(x)\ge\int_I \chi(s)\,ds=:c_1>0. \] Consecutive values have opposite signs, hence \[ |A_{t_{j+1}}f_\delta(x)-A_{t_j}f_\delta(x)| =a_{j+1}(x)+a_j(x)\ge2c_1. \] It follows that \[ V_q(Af_\delta;[1,2])(x)\gtrsim N^{1/q} \quad\text{on }B(0,c_0\delta). \] Thus \[ \|V_q(Af_\delta;[1,2])\|_p \gtrsim N^{1/q}\delta^{3/p}. \] 
Finally, choosing $N\sim \delta^{-1},$ it follows that 
\[ \frac{\|V_q(Af_\delta;[1,2])\|_p}{\|f_\delta\|_p} \gtrsim N^{1/q-1/p}\delta^{1/p}\gtrsim \delta^{2/p-1/q}.\] 
Taking $\delta\to 0,$ we obtain 
\[\frac1q\le\frac2p \iff q\ge\frac p2.\]    
\end{proof}

\begin{prop}
    Let \eqref{variationbound} holds for all $f\in L
    ^p$. Then $q$ must satisfy $q\geq3$.
\end{prop}
\begin{proof}
    The example below is adapted from Example $3.7$ of \cite{BOR}. Let $P_{j,n}=\Big\{(y_1,y_2,y_3)\in \mathbb{R}^3: |y_1|\leq 2^{-\frac{j}{2}-2}, |y_2|\leq 2^{-j-1}, |y_3-c_{j,n}|\leq 2^{-\frac{3j}{2}-3}\Big\}$, where $c_{j,n}=-n2^{-3j/2}$ for $n=1,2,\dots,N$ and $j>1$. Then $|P_{j,n}|\lesssim 2^{-3j}$. Choose $N\leq \lfloor{2^{\frac{3j}{2}-2}\rfloor}$. Consider 
    \begin{eqnarray*}
        &&f_j=\sum^N_{j=1}(-1)^n \chi_{P_{j,n}}, ~~\text{then}\\
        &&\Vert f_j\Vert_{L^p}\lesssim N^{1/p}|P_{j,n}|^{1/p}\lesssim 2^{\frac{3j}{2p}} 2^{-\frac{3j}{p}}= 2^{-\frac{3j}{2p}}.
    \end{eqnarray*}
    On the other hand, consider 
    \begin{eqnarray*}
        \Omega_j=\Big\{ (x_1,x_2,x_3)\in \mathbb{R}^3: |x_1|\leq 2^{-j/2}, |x_2|\leq 2^{-j}, 1\leq x_3\leq \frac{3}{2}\Big\}.
    \end{eqnarray*}
    Therefore, $|\Omega_j|\sim 2^{-3j/2}$.  Define $t_n(x)=|x_3-c_{j,n}|\in [1,2]$. Then, there exists a positive constant $a$ such that for $x\in \Omega_j$,
    \begin{eqnarray*}
        &&A_{t_{2\nu}}f_j(x)=\int_{I} f_j\Big(x-t_{2\nu}(x)\gamma(s)\Big) ds\gtrsim a2^{-j/2},\\
        &&A_{t_{2\nu+1}}f_j(x)=\int_{I} f_j\Big(x-t_{2\nu+1}(x)\gamma(s)\Big) ds\lesssim -a2^{-j/2}.
    \end{eqnarray*}
    Therefore, 
    \begin{eqnarray*}
        |V_q(\mathcal{A}f_j)(x)|\gtrsim N^{1/q}2^{-j/2}.
    \end{eqnarray*}
    For $N\sim {2^{3j/2-2}}$,
    \begin{eqnarray*}
        \Vert V_q(\mathcal{A}f_j;[1,2])\Vert_{L^p}\gtrsim N^{1/q}2^{-j/2} |\Omega_j|^{1/p}\gtrsim 2^{3j/2q-j/2-3j/2p}.
    \end{eqnarray*}
    Therefore, if the estimate \eqref{variationbound} holds,  then 
    \begin{eqnarray*}
        &&2^{3j/2q-j/2-3j/2p}\lesssim 2^{-3j/2p}\\
        &&\implies 2^{3j/2q}\lesssim 2^{j/2}\implies 3\leq q.
    \end{eqnarray*}
    This completes the proof.
\end{proof}

\section{Necessity of $q\geq \frac{6p}{p+4}$ in Theorem~\ref{mainresult}}   
\label{sec2}

\begin{prop}
\label{thm:ex:2}
 Let \eqref{variationbound} hold for all $f\in L^p.$ Then we necessarily have 
\begin{equation*}
\label{cond:2}
 q\geq \frac{6p}{p+4}.   
\end{equation*}   
\end{prop}

The proof of the above proposition is divided into several lemmas. Before that, we provide some known results and notations from \cite{BH}. Our example will be constructed for the curve $\gamma(s)=v+\left(s,\frac{s^2}{2},\frac{s^3}{6}\right),$ where $v=(1, -\frac{1}{2}, \frac{1}{6})$ and $\chi$ is a bump function centered at $s_0=0$. The translation is made in such a way that \begin{equation*}\label{eq:normalization}
\gamma^{(j)}(0)=e_j\ (1\le j\le3),\qquad
 a:=\gamma_3(0)>0.
\end{equation*}

We shall follow Beltran--Hickman \cite{BH} and make every notation precise for our prototypical example. Since $\gamma''(s)=(0,1,s)$, the equation
$\langle\gamma''(\theta(\xi)),\xi\rangle=0$ gives
\begin{equation*}\label{eq:theta}
\theta(\xi)=-\xi_2/\xi_3,\qquad \xi_3>0.
\end{equation*}
Solving simultaneously
$\langle\gamma'(s),\xi\rangle=\langle\gamma''(s),\xi\rangle=0$
with $\xi_3=1$ gives
\begin{equation*}\label{eq:Gamma}
\Gamma(\tau)=(\tau^2/2, \tau,1),\qquad
\theta(\Gamma(\tau))=-\tau.
\end{equation*}
Then the binormal cone is
$\mathcal C=\{r\Gamma(\tau):r>0,\ |\tau|\le c_0\}$.
Define
\begin{equation*}\label{eq:phi}
\begin{split}
\phi(\xi)&=\langle\gamma(\theta(\xi)),\xi\rangle
=v\cdot\xi-\frac{\xi_1\xi_2}{\xi_3}+\frac{\xi_2^3}{3\xi_3^2},\\
u_3(\xi)&=\langle\gamma'''(\theta(\xi)),\xi\rangle=\xi_3.
\end{split}
\end{equation*}
On a fixed cone where $\xi_3\ge c|\xi|>0$, the function $\phi$ is smooth and homogeneous of degree one, leading to 
\begin{equation}\label{eq:phase-derivatives}
|\partial_\xi^\beta\phi(\xi)|\lesssim_\beta |\xi|^{1-|\beta|}.
\end{equation}
Differentiating the definition gives
\[
\nabla\phi(\xi)=\gamma(\theta(\xi))+
\langle\gamma'(\theta(\xi)),\xi\rangle\nabla\theta(\xi).
\]
The second term vanishes on $\mathcal C$, so
\begin{equation}\label{eq:gradient}
\nabla\phi(\lambda\Gamma(\tau))=\gamma(-\tau).
\end{equation}
Moreover,
\begin{equation*}\label{eq:temporal-phase}
\phi(\Gamma(\tau))=\frac{(1-\tau)^3}{6}.
\end{equation*}
Choose $0<c_0\le1/2$ sufficiently small that
\begin{equation}\label{imp:est:1}
\chi(-\tau)\ge c_\chi>0,\qquad
\left|\frac{\phi(\Gamma(\tau))}{a}-1\right|\le\frac1{10},
\quad (|\tau|\le c_0).
\end{equation}

Set $\delta=\lambda^{-1/3}$ and
\begin{equation*}\label{eq:indices}
\begin{split}
\NN(\lambda)&=\Z\cap[-c_0\lambda^{1/3},c_0\lambda^{1/3}],\\
\xi^\nu&=\lambda\Gamma(\nu\delta)
=\left(\frac{\nu^2}{2}\lambda^{1/3},\nu\lambda^{2/3},\lambda\right),
\qquad s_\nu=-\nu\delta.
\end{split}
\end{equation*}
Thus $\#\NN(\lambda)\sim \lambda^{1/3}$.
Choose $\eta\in C_c^\infty(B(0,1))$ be a non-negative function, which is equal to $1$ on $B(0,1/2)$,
and $\eta_+\in C_c^\infty(B(0,3/2))$ equal to $1$ on $B(0,1)$.
For a small fixed $\rho>0$, put $h=\rho\lambda^{1/3}$ and define
\begin{equation}\label{eq:g}
\widehat g_\nu(\xi)=\eta(\frac{(\xi-\xi^\nu)}{h}),\qquad
\widehat g_{\nu,+}(\xi)=\eta_+(\frac{(\xi-\xi^\nu)}{h}).
\end{equation}
Then $\widehat g_{\nu,+}\widehat g_\nu=\widehat g_\nu$, and
\begin{equation}\label{eq:gnorm}
\|g_\nu\|_2^2=(2\pi)^{-3}\rho^3\lambda\int|\eta(\zeta)|^2\,d\zeta
\sim c_{\rho}\lambda.
\end{equation}
Note that $|\xi^\nu-\xi^{\nu'}|\ge|\nu-\nu'|\lambda^{2/3}$, the balls $B(\xi^\nu, \rho \lambda^{1/3})$ are well-separated. 

Fix $t_*=3/2$ and $0<c_1<1/4$, and set
\begin{equation*}\label{eq:times}
\MM(\lambda)=\{m\in\Z:|m|\le c_1\lambda^{1/3}\},\qquad t_m=t_*+m\delta.
\end{equation*}
Then $\#\MM(\lambda)\sim \lambda^{1/3}$ and all $t_m$ lie in a fixed
compact subinterval of $(1,2)$. Define
\begin{equation*}\label{eq:fpacket}
\widehat f_{m,\nu}(\xi)=\lambda^{1/3}e^{it_m\phi(\xi)}\widehat g_\nu(\xi),
\qquad F_\lambda^\omega=\sum_{m\in\MM(\lambda)}\sum_{\nu\in\NN(\lambda)}
\epsilon_{m,\nu}(\omega)f_{m,\nu},\quad \omega\in [0, 1],
\end{equation*}
and $\{\epsilon_{m,\nu}\}$ sequence of Rademacher functions. 
Define  \begin{align*}
c_{t}(\xi)= \alpha \lambda^{1/3} \chi(\theta(\xi)) (t u_3(\xi))^{-1/3}\qquad c_{t,\nu}:=c_{t}(\xi^\nu)=\alpha \chi(\theta(\xi^\nu))\frac{\lambda^{1/3}}{(t u_3(\xi^\nu))^{1/3}}
=\alpha\chi(s_\nu)t^{-1/3}.
\label{eq:ctnu}
\end{align*}
The following lemma is established in \cite{BH}. We state it here for as it will play a crucial role in our construction.
\begin{lemma} \cite[Lemma 3.2]{BH}\label{asy-lemma} 
Let $1 \leq t \leq 2$. If $c_0$, $\rho > 0$ are chosen sufficiently small, then for all $\lambda \geq 1$ and $\nu\in \NN(\lambda)$ the identity
\begin{equation*}
    \widehat{\mu}_{t}(\xi)= e^{-i t\phi(\xi)} m_t(\xi)
\end{equation*} 
holds on $\supp \widehat{g}_{+,\nu}$ where
\begin{enumerate}[i)]
    \item The function $m_t$ satisfies the asymptotics
\begin{align}
\label{asymp:1}
|\lambda^{1/3} m_t(\xi)-c_{t}(\xi)|\lesssim_{\chi} (\rho + \lambda^{-1/3}), \quad \text{for}\qquad \xi \in \supp \widehat{g}_{\nu,+}.
\end{align}
\smallskip
    \item The function $m_t$ satisfies the derivative bounds
    \begin{equation}\label{eq:der}
    |\partial^\alpha_\xi m_t(\xi)| \lesssim \lambda^{-1/3 - |\alpha|/3}  \quad \text{ for all $\alpha \in \N_0^n$}
    \end{equation}
    and $\xi \in \supp \widehat{g}_{\nu,+}$.
 \item Moreover, it follow from \cite{BH}(just above equation 3.10 in \cite{BH}) that
\begin{align*}
|c_{t}(\xi)-c_{t}(\xi^\nu)| \lesssim \lambda^{-2/3},\quad \xi \in \supp \widehat{g}_{\nu,+}.    
\end{align*}
Consequently, 
\begin{align}
\lambda^{1/3} m_{t}(\xi)=c_{t, \nu}+E_{t,\nu, \lambda}(\xi),\quad \xi \in \supp \widehat{g}_{\nu,+},
\label{use:est}
\end{align}
where $|E_{t,\nu, \lambda}(\xi)|=O(\rho + \lambda^{-1/3})$ for $\xi \in \supp \widehat{g}_{\nu,+}.$
\end{enumerate}
\end{lemma}
\begin{lemma}\label{sep:centre}
Let $z_{m,\nu}=-t_m\gamma(s_\nu)$.
Then for distinct pairs $(m,\nu)$ and $(m',\nu')$, we have 
\[
|z_{m,\nu}-z_{m',\nu'}|
\gtrsim \lambda^{-1/3} \bigl(|m-m'|+|\nu-\nu'|\bigr)\ge c\lambda^{-1/3}.
\]
\end{lemma}
\begin{proof}
The proof follows from obvious considerations and the inequality \eqref{sep:ex:1}. 
\end{proof}
\begin{lemma}
For every integer $N\ge1$,
\begin{equation}\label{eq:packetdecay}
|f_{m,\nu}(x)|\le C_N\lambda^{4/3}
\left(1+\lambda^{1/3}|x+t_m\gamma(s_\nu)|\right)^{-N}.
\end{equation}
Moreover, for $1\leq p\leq\infty$ we have
\begin{equation}\label{lp:norm}
\|F_\lambda^\omega\|_p\lesssim\lambda^{4/3-1/(3p)},
\end{equation}
 uniformly over all choices of signs.
 \end{lemma}

\begin{proof}
Using inversion formula we get 
\begin{align*}
\nonumber f_{m, \nu}(x)&\sim \int e^{i\, x\cdot \xi} \widehat f_{m\nu}(\xi)\, d\xi\\
&= \lambda^{1/3}\int e^{ix.\xi} e^{it_{m}\phi(\xi)} \eta\big(\frac{\xi-\xi^{\nu}}{\rho \lambda^{1/3}}\big)~d\xi\nonumber\\
&= \lambda^{1/3}\int e^{ix.(\xi^\nu+\xi \lambda^{1/3}\rho)} e^{it_{m}\phi(\xi^\nu+\lambda^{1/3}\rho\xi)} \eta(\xi) \lambda \rho^3\, d\xi
\end{align*}
Denote $h:=\lambda^{1/3}\rho,$ by Taylor's expansion we obtain
\[
\phi(\xi^\nu+h\xi)=\phi(\xi^\nu)+h\nabla \phi(s_\nu).\xi+ R_\nu(\xi).
\]
By \eqref{eq:gradient},$\nabla \phi(\xi^\nu)=\gamma(s_{\nu}),$ hence 
\[
\phi(\xi^\nu+h\xi)=\phi(\xi^\nu)+h\gamma(s_{\nu}).\xi+ R_\nu(\xi).
\]
Also, $\phi$ is homogeneous of degree 1, hence near $\xi^\nu,$ $|\partial^\alpha_{\xi}\phi(\tilde{\xi})|\lesssim |\lambda|^{1-|\alpha|}.$ On $|\xi|\le1$, the remainder and its first derivatives are $O(h^2/\lambda)$. For $|\alpha|\ge2$, \eqref{eq:phase-derivatives} gives
\[
|\partial_\zeta^\alpha R_\nu(\zeta)|\lesssim_\alpha h^{|\alpha|}\lambda^{1-|\alpha|}.
\]
All these quantities are bounded uniformly as $\lambda\to\infty$. Therefore, 
\begin{align*}
f_{m\nu}(x)&\sim \lambda^{1/3}\int e^{ix.(\xi^\nu+\xi \lambda^{1/3}\rho)} e^{it_{m}\phi(\xi^\nu+\lambda^{1/3}\rho\xi)} \eta(\xi) \lambda \rho^3\, d\xi\\
&=\lambda^{1/3} h^3
 e^{i(x\cdot\xi^\nu+t_m\phi(\xi^\nu))}\int e^{ih(x+t_m\gamma(s_\nu))\cdot\xi}
 e^{it_mR_\nu(\xi)}\eta(\xi)\,d\xi.
\end{align*}

The above consideration together with integration by parts now conclude that 
$$|f_{m,\nu}(x)|\le C_N\lambda^{4/3}
\left(1+\lambda^{1/3}|x+t_m\gamma(s_\nu)|\right)^{-N}.$$
For the second part, 
choose an integer $N>3$. Integrating \eqref{eq:packetdecay} gives
\[
\|F_\lambda^\omega\|_1
\lesssim\lambda^{4/3}\#\MM(\lambda)\#\NN(\lambda)\lambda^{-1}
\int_{\R^3}(1+|y|)^{-N}\,dy
\lesssim \lambda^{4/3} \lambda^{1/3}\, \lambda^{1/3} \lambda^{-1}\sim \lambda.
\]

For the $L^\infty$ estimate, set
\[
H(x)=\sum_{m,\nu}(1+\lambda^{1/3}|x-z_{m,\nu}|)^{-N}.
\]
For each $\ell\ge0$, the annulus
$2^\ell \lambda^{-1/3} \le|x-z_{m,\nu}|<2^{\ell+1}\lambda^{-1/3}$
contains at most $C 2^{3\ell}$ centres, each contributing at most
$2^{-N\ell}$. Therefore
\[
H(x)\lesssim1+\sum_{\ell\ge0}2^{-(N-3)\ell}\lesssim1.
\]
Since $|F_\lambda^\omega(x)|\lesssim\lambda^{4/3}H(x)$, this proves
$\|F_\lambda^\omega\|_\infty\lesssim\lambda^{4/3}$.
For $1\le p<\infty$, by interpolation we derive
\[
\|F_\lambda^\omega\|_p^p
\le\|F_\lambda^\omega\|_\infty^{p-1}\|F_\lambda^\omega\|_1
\lesssim\lambda^{(4/3)(p-1)}\lambda
=\lambda^{4p/3-1/3}.
\]
This completes the proof.
\end{proof}

By definition,
\begin{equation*}\label{eq:outputfourier}
\widehat{A_tf_{m,\nu}}(\xi)
=e^{-i(t-t_m)\phi(\xi)}\lambda^{1/3}m_t(\xi)\widehat g_\nu(\xi).
\end{equation*}
Fix $c_2>0$ and let $J_m=[t_m-c_2\lambda^{-1/3},t_m+c_2\lambda^{-1/3}]$.
For sufficiently large $\lambda$, these intervals lie in a fixed compact
subinterval of $(1,2)$.
\begin{lemma}\label{lem:output}
For $t\in J_m$ and every integer $N\ge1$,
\begin{equation}\label{eq:outputdecay}
|A_tf_{m,\nu}(x)|\le C_N\lambda(1+\lambda^{1/3}|x|)^{-N}.
\end{equation}
\end{lemma}

\begin{proof}
Again write $h=\rho\lambda^{1/3}$ and $\xi=\xi^\nu+h\zeta$. Using Fourier inversion we obtain
\[
A_tf_{m,\nu}(x)\sim h^3
 e^{i x\cdot\xi^\nu-i(t-t_m)\phi(\xi^\nu)}
\int e^{ihx\cdot\zeta}B_{m,\nu,t}(\zeta)\,d\zeta,
\]
where
\begin{align*}
B_{m,\nu,t}(\zeta)={}&e^{-i(t-t_m)[\phi(\xi^\nu+h\zeta)-\phi(\xi^\nu)]}\lambda^{1/3}m_t(\xi^\nu+h\zeta)\eta(\zeta).
\end{align*}
The multiplier factor has uniformly bounded derivatives by \eqref{eq:der}. For $r\ge1$, a derivative of order $r$ of the phase in brackets, multiplied by $t-t_m$, is bounded by
\[
C_r \lambda^{-1/3} h^r\lambda^{1-r}\lesssim_r1.
\]
Thus every derivative of $B_{m,\nu,t}$ is uniformly bounded, and its support is fixed. Integration by parts yields $|A_{t}(f_{m\nu})(x)|\lesssim C_N\, h^3(1+h|x|)^{-N}$. Since $h=\rho\lambda^{1/3}$, this gives \eqref{eq:outputdecay}.
\end{proof}

We shall now proceed with computing variational seminorm. Recall $a=\gamma_{3}(0).$ Let us define
$$\tau_j=1+j\frac{\pi}{a\lambda}\quad, 0\le j\le K, K\sim \lfloor\lambda \rfloor,$$
and $$D_jf=A_{\tau_{j+1}}f-A_{\tau_j}f$$ for $0\le j<K$. Moreover, recall from \eqref{imp:est:1} that  $\left|\frac{\phi(\xi^\nu)}{a\lambda}-1\right|\le\frac1{10},$ which by continuity also implies 
\begin{align}
\label{imp:est:2}
\bigg|\frac{\phi(\xi)}{a\lambda}-1\bigg|\lesssim \frac{1}{10},\quad \xi \in \supp \widehat{g}_{\nu,+}.   
\end{align}
Now, from \eqref{imp:est:2} it follows that $$|e^{-i\frac{\phi(\xi)}{a\lambda}}-1|=2|\sin\left(\frac{\pi\phi(\xi)}{2a\lambda}\right)|
\ge \kappa>0,$$
where this $\kappa$ is a fixed constant independent of $\nu, \xi\in \supp \widehat{g}_{\nu,+}.$

For $m\in \MM(\lambda),$ we denote the set 
\[
\mathcal J_m=\{j:0\le j<K,\ [\tau_j,\tau_{j+1}]\subset J_m\}.
\]
Then $\# \mathcal J_m\sim \lambda^{2/3}.$
\begin{lemma}
Using the notation introduced above, we derive the following estimate
 \begin{equation*}
 \label{lower:estimate}
\| D_j f_{m,\nu} \|_{2}\gtrsim \kappa \lambda^{1/2}.   
 \end{equation*}
 Additionally, 
 let $m\in \MM(\lambda),$ and $j\in \mathcal J_m.$ There exists large $R>0,$ independent of indices and $\lambda$, such that 
\begin{align*}
\label{D:local}
\int_{E_\lambda}|D_jf_{m, \nu}(x)|^2dx\gtrsim \frac{\kappa^2}{2} \lambda,    
\end{align*}
where $E_{\lambda}=B(0, R\lambda^{-1/3}).$
\end{lemma}

\begin{proof}
Applying the Fourier transform to $D_j f_{m,\nu}$, we obtain  
\begin{equation*}
\widehat{D_j f_{m,\nu}}(\xi)
=
e^{-i(\tau_j-t_m)\phi(\xi)}
Q_{j,\nu}(\xi)\widehat{g}_\nu(\xi),
\end{equation*}
where
\begin{equation*}
Q_{j,\nu}(\xi)
=
e^{-i\frac{\phi(\xi)}{a\lambda}\pi}
\lambda^{1/3}m_{\tau_{j+1}}(\xi)
-
\lambda^{1/3}m_{\tau_j}(\xi).
\end{equation*}
Using \eqref{use:est} at both times and adding and subtracting
$e^{-i\frac{\phi(\xi)}{a\lambda}}c_{\tau_j,\nu}$, we obtain
\begin{align*}
Q_{j,\nu}(\xi)&= c_{\tau_j,\nu}
\bigl(e^{-i\frac{\phi(\xi)}{a\lambda}}-1\bigr)
+ e^{-i\frac{\phi(\xi)}{a\lambda}}
\bigl(c_{\tau_{j+1},\nu}-c_{\tau_j,\nu}\bigr)
\\
&+
\left(e^{-i\frac{\phi(\xi)}{a\lambda}}
E_{\tau_{j+1},\nu,\lambda}(\xi)-
E_{\tau_j,\nu,\lambda}(\xi)\right)\nonumber\\
&=I_1(\xi)+I_2(\xi)+I_3(\xi),
\end{align*}

where, as a consequence of \eqref{imp:est:2}, we further estimate 
\begin{align}
\label{i1}
|I_1|=|c_{\tau_j,\nu}
\bigl(e^{-i\frac{\phi(\xi)}{a\lambda}}-1\bigr)|=|c_{\tau_j,\nu}| \big|e^{-i\frac{\phi(\xi)}{a\lambda}}-1\big|\gtrsim |\chi(s_{\nu})| |\tau_{j}^{-1/3}|\gtrsim \kappa,    
\end{align}
where this $\kappa>0$ is independent of $\nu, \xi\in \supp \widehat{g_{\nu,+}}.$

Since $|\partial_tc_{t,\nu}|=|-(\alpha/3)\chi(s_\nu)t^{-4/3}|\lesssim 1,$ we obtain $|c_{\tau_{j+1},\nu}-c_{\tau_j,\nu}|\lesssim |\tau_{j+1}-\tau_j|\lesssim \frac{1}{\lambda},$ and as a consequence, we derive 
\begin{align}
\label{est:i2}
|I_2|\lesssim \frac{1}{\lambda}.   
\end{align}
Moreover, \eqref{asymp:1} gives
\[
\bigl|E_{t,\nu,\lambda}(\xi)\bigr|
\lesssim \rho+\lambda^{-1/3}
\]
uniformly for $t\in[1,2]$ and
$\xi\in\operatorname{supp}\widehat{g}_\nu$.
Consequently, using \eqref{i1} and \eqref{est:i2}, we derive
\begin{align*}
|Q_{j,\nu}(\xi)|\geq|I_1|-|I_2|-|I_3|\gtrsim \frac{\kappa}{2},
\end{align*}
by taking $\lambda$ very large and $\rho$ sufficiently small. Then $|Q_{j,\nu}|\ge \frac{\kappa}{2}>0$ on $\supp\widehat g_\nu$. Using Plancherel theory and \eqref{eq:gnorm} we get 
\begin{equation}\label{eq:Dglobal}
\|{D_jf_{m,\nu}}\|_2^2\gtrsim \kappa^2 \|g_\nu\|_2^2\gtrsim \kappa^2 \lambda.
\end{equation}
This completes the proof of the first estimate.

Let $j\in\mathcal J_m$. Since Lemma~\ref{lem:output} holds at both endpoints. It yields
\[
\int_{|x|>R\lambda^{-1/3}}|D_jf_{m,\nu}(x)|^2dx
\le C_N\lambda\int_{|y|>R}(1+|y|)^{-2N}dy.
\]
Choose $N>3/2$ and a sufficiently large fixed $R$ so that the right side is at most $\kappa^2\lambda/2$. With
\begin{equation*}\label{eq:region}
E_\lambda=B(0,R\lambda^{-1/3}),\qquad |E_\lambda|\sim\lambda^{-1},
\end{equation*}
we conclude from \eqref{eq:Dglobal} that
\begin{equation}\label{eq:Dlocal}
\int_{E_\lambda}|D_jf_{m,\nu}(x)|^2dx\ge \frac{\kappa^2}{2}\lambda,
\qquad j\in\mathcal J_m.
\end{equation}
\end{proof}

The consequence of the above lemma implies that most of the mass of $|D_j f_{m, \nu}|$ is localised on the ball $E_{\lambda}=B(0, R\lambda^{-1/3}).$

\begin{lemma}
Let $\lambda\gg 1.$ Then there exists a choice of sign $\omega_{\lambda}=\{\epsilon_{m, \nu}\},$ depending on $\lambda,$ such that 
$$\int_{E_\lambda}\sum_{j=0}^{K-1}|D_{j} F_{\lambda}^{\omega}(x)|^2\, dx\gtrsim \lambda^{7/3},$$
where $F^{\omega}_{\lambda}=\sum\limits_{m\in \MM(\lambda)}\sum\limits_{\nu\in \NN(\lambda)} \epsilon_{m, \nu} f_{m, \nu}.$
\end{lemma}

\begin{proof}
Consider 
\begin{align*}
\mathbb E|D_jF_\lambda^\omega(x)|^2
&=\sum_{m,\nu}\sum_{m',\nu'}
\mathbb E(\epsilon_{m,\nu}\epsilon_{m',\nu'})
D_jf_{m,\nu}(x)\overline{D_jf_{m',\nu'}(x)}\\
&=\sum_{m,\nu}|D_jf_{m,\nu}(x)|^2.
\end{align*}
Define
\[
I_\lambda(\omega)=\int_{E_\lambda}
\sum_{j=0}^{K-1}|D_jF_\lambda^\omega(x)|^2\,dx.
\]
Using \eqref{eq:Dlocal} we obtain
\begin{align}
\mathbb E I_\lambda
&=\sum_{m\in\MM(\lambda)}\sum_{\nu\in\NN(\lambda)}\sum_{j=0}^{K-1}
\int_{E_\lambda}|D_jf_{m,\nu}(x)|^2\,dx\nonumber\\
&\ge\sum_{m\in\MM(\lambda)}\sum_{\nu\in\NN(\lambda)}\sum_{j\in\mathcal J_m}
\int_{E_\lambda}|D_jf_{m,\nu}(x)|^2\,dx\nonumber\\
&\gtrsim\lambda^{1/3}\lambda^{1/3}\lambda^{2/3}\lambda
=\lambda^{7/3}.\label{expec:lower}
\end{align}
Consequently, there exists $\omega_\lambda$ such that, for 
$F_\lambda=F_\lambda^{\omega_\lambda}$ we have
\begin{equation*}\label{eq:selected-energy}
\int_{E_\lambda}\sum_{j=0}^{K-1}|D_jF_\lambda(x)|^2\,dx
\gtrsim\lambda^{7/3}.
\end{equation*}
\end{proof}

Now we are at a position to complete the proof of Proposition~\ref{thm:ex:2}. 

\begin{proof}[Proof of Proposition~\ref{thm:ex:2}]
Recall $K\sim \lfloor\lambda \rfloor,$ and $q\geq 2.$ Using H\"older's inequality, we observe that 
\begin{align*}
\nonumber V_q(AF_\lambda;[1,2])(x)&\geq \left(\sum_{j=0}^{K-1}|A_{\tau_{j+1}}F_\lambda(x)-A_{\tau_{j}}F_\lambda(x)|^q\right)^{1/q}\\
&\geq \lambda^{\frac{1}{q}-\frac{1}{2}}\left(\sum_{j=0}^{K-1}|A_{\tau_{j+1}}F_\lambda(x)-A_{\tau_{j}}F_\lambda(x)|^2\right)^{1/2}\sim \lambda^{\frac{1}{q}-\frac{1}{2}} \left(\sum_{j=0}^{K-1}|D_j F_{\lambda}(x)|^2\right)^{1/2}.
\end{align*}

Using \eqref{expec:lower} we further conclude that 
\begin{align*}\label{lower:L2}\|V_q(AF_\lambda;[1,2])\|_{L^2(E_\lambda)}^2&\gtrsim \lambda^{\frac{2}{q}-1}\lambda^{7/3}\sim\lambda^{2/q-1+7/3}=\lambda^{\frac{2}{q}+\frac{4}{3}}.\end{align*}

Again using H\"older's inequality and $|E_\lambda|\sim \lambda^{-1},$ we conclude 
\begin{equation}
\|V_q(AF_\lambda;[1,2])\|_{L^p(E_\lambda)}
\ge
|E_\lambda|^{1/p-1/2}\|V_q(AF_\lambda;[1,2])\|_{L^2(E_\lambda)}\gtrsim \lambda^{1/2-1/p} \lambda^{1/q+2/3}. \label{final:lower:bound}     
\end{equation}

From \eqref{lp:norm} and \eqref{final:lower:bound}, we finally derive 

\begin{align*} \frac{\|V_q(AF_\lambda;[1,2])\|_p}{\|F_\lambda\|_p} &\gtrsim \lambda^{1/q+7/6-1/p-(4/3-1/(3p))}\gtrsim \lambda^{1/q-1/6-2/(3p)}. \end{align*} 

Taking $\lambda\to \infty,$ we obtain $\frac1q-\frac16-\frac{2}{3p}\leq 0\iff q\ge\frac{6p}{p+4}.$
This completes the proof.
    
\end{proof}

\section{Proof of sufficiency part of Theorem \ref{mainresult}}
\label{suff}
\subsection{Preparatory results}
Since the long-variation already satisfies the desired estimate, the proof of Theorem~\ref{mainresult} requires us to show that the short variation
\begin{align*}
S_{q}(\mathcal Af)(x)
=
\Big(
\sum_{j\in\mathbb Z}
\big[V_{q,j}(\mathcal Af)(x)\big]^q
\Big)^{1/q}\end{align*}
is bounded from $L^p(\R^3)$ to itself if $q>\max\{3, \frac{6p}{p+4}, \frac{p}{2}\}.$

To start with, we will dyadically decompose $S_{3}(\mathcal Af).$ By rescaling
\begin{align*}
S_3(\mathcal A f)(x)
=
\Big(
\sum_{j\in\mathbb Z}
V_{3,j}\bigl(A_{t}f(x):t\in[1,2]\bigr)^3
\Big)^{1/3}=\left(
\sum_{j\in\mathbb Z}
V_{3}\big(A_{2^jt}f(x):t\in[1,2]\big)^3
\right)^{1/3}.
\end{align*}
Recall the Littlewood-Paley projections $\widehat{P_{k}f}({\xi})=\widehat{\beta_{k}(D)f}(\xi)=\beta_{k}(\xi) \widehat{f}(\xi),$ where $\beta_{k}(\xi):=\beta(|\xi|/2^{k})$ for $k\geq1$ is defined earlier, and $\beta_0(r)=\zeta (r)$. Also, $\widetilde P_k$ will denote Littlewood-Paley projections on slightly larger
annulus such that $\widetilde P_k P_k=P_k.$

For each fixed \(j\), writing
$f=\sum_{k\geq 0}P_{k-j}f$, we get
\begin{align*}
V_3\big(A_{2^jt}f:t\in[1,2]\big)
\leq
\sum_{k\geq0}
V_3\big(A_{2^jt}P_{k-j}f:t\in[1,2]\big).    
\end{align*}
Consequently, using Minkowski's inequality, we obtain
\begin{align*}
S_3(\mathcal A f)(x)
&\leq
\big(
\sum_j
\big[
\sum_{k\geq0}
V_3\big(A_{2^jt}P_{k-j}f(x):t\in[1,2]\big)
\big]^3
\big)^{1/3}\nonumber\\
&\leq
\sum_{k\geq0}
\big(
\sum_j
V_3\big(A_{2^jt}P_{k-j}f(x):t\in[1,2]\big)^3
\big)^{1/3}=\sum_{k\geq0} S_{3}(\mathcal A^k f),
\end{align*}
where $$S_{3}(\mathcal A^k f)(x):=\Big(
\sum_j
V_3\big(A_{2^jt}P_{k-j}f(x):t\in[1,2]\big)^3
\Big)^{1/3}.$$

For the \textit{low frequency} $k=0,$ we obtain even a better estimate. Using fundamental theorem of calculus and H\"older's inequality we obtain 
\begin{align*}
S_{2}(\mathcal{A}^0f)(x)\lesssim \left(\int_{0}^{\infty}\big|t \partial_{t} (\A^{0}_{t}f) \big|^2\, \frac{dt}{t}\right)^{1/2}.
\end{align*}
Note that, $t \widehat{\partial_{t} (\A^{0}_{t}f)}(\xi)=\widehat{f}(\xi)\widehat{\Psi(t\xi)},$ where $\widehat{\Psi(\xi)}=\langle \xi\cdot \nabla_{\xi} (\widehat{\mu}(\xi)\beta_0(|\xi|)) \rangle.$ Since $\beta_0$ is a compactly supported smooth function and $\hat{\mu}$ is a smooth function in $\xi,$ therefore, $\Psi$ is a Schwartz class function with $\widehat{\Psi}(0)=0,$ thus appealing to the vector-valued Calder\'on-Zygmund theory, see \cite[Chapter I]{Steinfat},  we derive
\begin{align*}
\|S_{2}(\mathcal{A}^0f)\|_{L^p(\R^3)}\lesssim \big\| \left(\int_{0}^{\infty}\big|t \partial_{t} (\A^{0}_{t}f) \big|^2\, \frac{dt}{t}\right)^{1/2} \big\|_{L^p(\R^3)}&\sim \big\|\big(\int_{0}^{\infty} |f\ast (t^{-3} \Psi(t^{-1}\cdot))|^2\frac{dt}{t}\big)^{1/2}\big\|_{L^p(\R^3)}\\
&\lesssim \|f\|_{L^p(\R^3)}.   
\end{align*}
Hence, for the \textit{low frequency} we obtain a better estimate. Next we obtain estimates for the dyadically localized $S_{3}(\mathcal A^k f)$ for $k\geq 1.$ 

\begin{theorem} Let $k\geq1$ and $3<p\leq 4$. Then for any $\epsilon>0$,
$$\|S_3(\mathcal{A}^k f)\|_{L^p} \leq C 2^{k\epsilon} \|f\|_{L^p}.$$
\label{no-decay}
\end{theorem}

The proof of the above Theorem requires some ideas from Proposition~7.2 \cite{BOR} with important modifications pertaining our situation. To begin with, by the Sobolev embedding and Young's inequality we obtain
\begin{align}
\label{embed}
V_3\big(A_{2^jt}P_{k-j}f(x):t\in[1,2]\big)
\lesssim
2^{k/3}
\Big(\int_1^2|A_{2^jt}P_{k-j}f(x)|^3\,dt\Big)^{1/3}+
2^{-2k/3}
\Big(
\int_1^2
\big|\partial_t[A_{2^jt}P_{k-j}f(x)]\big|^3\,dt
\Big)^{1/3},
\end{align}

which further implies the estimate

\begin{align}
\label{ptwise:dyadic:sv}
S_3(\A^kf)(x)
\lesssim
2^{k/3}
\Big(
\sum_j\int_1^2
|A_{2^jt}P_{k-j}f(x)|^3\,dt
\Big)^{1/3}+
2^{-2k/3}
\Big(
\sum_j\int_1^2
\left|\partial_t[A_{2^jt}P_{k-j}f(x)]\right|^3\,dt
\Big)^{1/3}.
\end{align}

The following lemma combines individual dyadic pieces. 

\begin{prop}
\label{global:prop}
Fix $k\geq 1.$ For \(p>3\), we have
\begin{align*}
&\big\|
\big(
\sum_j\int_1^2
|A_{2^jt}P_{k-j}f|^3\,dt
\big)^{1/3}
\big\|_{L^p(\R^3)}+2^{-k}
\big\|
\big(
\sum_j\int_1^2
\left|\partial_t[A_{2^jt}P_{k-j}f]\right|^3\,dt
\big)^{1/3}\big\|_{L^p(\R^3)}\lesssim c_{p,\epsilon}
2^{-k(\sigma(p)-\epsilon)}\big\|f\big\|_{L^p(\R^3)}.
\end{align*}
\end{prop}

\begin{proof}
We denote
\begin{align*}
G(x)=\big(
\sum_j\int_1^2
|A_{2^jt}P_{k-j}f(x)|^3\,dt
\big)^{1/3},
\end{align*}
and the proof uses Fefferman--Stein maximal function to estimate $\|G\|_{L^p(\R^3)}.$  Here on wards, we simply write $\mathcal{A}^{k}_{j}f(x,t)=A_{2^jt}P_{k-j}f(x).$ The sharp maximal function of $G$ is defined as
$$G^{\#}(x)=\sup_{Q\ni x}\frac{1}{|Q|}\int_{Q}|G(y)-G_{Q}|\, dy,\quad \text{where}\quad G_Q=\frac{1}{|Q|}\int_Q G,$$
and the supremum is over all cubes $Q,$ with sides parallel to the coordinate axes, containing $x.$ For $L\in\Z,$  $\mathcal Q_L(x)$ will denote the collection of all cubes $Q,$ with sides parallel to the coordinate axes, containing $x$ such that $2^{L}\leq \ell(Q)< 2^{L+1}.$ This yields
\begin{equation}
 G^\#(x)\le 2\mathcal Gf(x)+\sum_{n\ge1}U_nf(x),
 \label{eq:split}
\end{equation}
where
\begin{align*}
 \mathcal Gf(x)&=\sup_{L\in\Z}\sup_{Q\in\mathcal Q_L(x)}
\frac{1}{|Q|}\int_{Q} \big(\sum_{j\le L}\|\mathcal{A}^k_j f(y,\cdot)\|_{L^3_t}^3\big)^{1/3}dy,\\
 U_nf(x)&=\sup_{L\in\Z}\sup_{Q\in\mathcal Q_L(x)}
 \frac{1}{|Q|^2}\int_Q\int_Q
 \|\mathcal A^k_{L+n}f(y,\cdot)-\mathcal{A}^k_{L+n}f(w,\cdot)\|_{L^3_t}\,dw\,dy .
\end{align*}
\medskip

\noindent\textbf{Estimates for $\mathcal {G}$:} Let $\mathfrak{M}f$ denotes the Hardy-Littlewood maximal function of $f.$ Trivially, using local smoothing estimates from Theorem~\ref{LS:thm} we obtain

\begin{align}
\nonumber \|\mathcal Gf\|_{L^3(\R^3)}&\lesssim \|\mathfrak{M}(\big(\sum_{j\in \Z}\|\mathcal{A}^k_j f(\cdot,\cdot)\|_{L^3_t}^3\big)^{1/3})\|_{L^3(\R^3)}\\
\nonumber&\lesssim \big(\sum_j\|\A^k_j(\widetilde{P_{k-j}} f\|_{L^3_x(L^3_t)}^3\big)^{1/3}\\
&\lesssim 2^{-k(\sigma(3)-\epsilon)} \big(\sum_j\|\widetilde P_{k-j}f\|_{L^3(\R^3)}^3\big)^{1/3}\lesssim 2^{-k(\sigma(3)-\epsilon)} \|f\|_{L^3(\R^3)},\label{l3est}
\end{align}
where the last inequality follows from Littlewood-Paley theorem.

To obtain an $L^\infty$ estimate for $\mathcal{G}f,$ fix $x,$ and a cube $Q\in \mathcal Q_L(x)$ for some fixed $L\in \Z.$ Then by triangle inequality, H\"older's inequality and $\ell^1\subset \ell^3$ we derive
\begin{align}
\nonumber&\frac{1}{|Q|}\int_{Q} \big(\sum_{j\le L}\|\mathcal{A}^k_j f(y,\cdot)\|_{L^3_t}^3\big)^{1/3}dy\\
&\lesssim \bigg(\frac{1}{|Q|} \int_{Q}\sum_{j\le L}\|\mathcal{A}^k_j (\chi_{B_x^L} f)(y,\cdot)\|_{L^3_t}^3\, dy\bigg)^{1/3}+\frac{1}{|Q|}\int_{Q}\sum_{j\le L}\|\mathcal{A}^k_j (\chi_{(B_x^L)^c} f)(y,\cdot)\|_{L^3_t}\, dy,\label{loc-glo} 
\end{align}
where $B_x^L=B(x,C_02^L)$ with $C_0$ sufficiently large. Using Theorem~\ref{LS:thm} the local part is bounded by 
\begin{align}
2^{-3L/3}2^{-k(\sigma(3)-\epsilon)}
 \left(\sum_{j\le L}\|\widetilde P_{k-j}(f\chi_{B_x^L})\|_{L^3(\R^3)}^3
 \right)^{1/3}
 \lesssim 2^{-k(\sigma(3)-\epsilon)} 2^{-3L/3} \|(f\chi_{B_x^L})\|_{L^3(\R^3)} \lesssim 2^{-k(\sigma(3)-\epsilon)} \|f\|_{L^\infty(\R^3)}.\label{loc} \end{align}

For the global part, write $\A^{k}_{j}f(x, t)=K^k_{j}\ast f(x),$ where the kernel $K^k_{j}(x)=2^{-3j} K^{k, t}(2^{-j}x)$ and $K^{k, t}$ is the convolution kernel for $A_{t}P_{k}.$ Recalling the estimate $|K^{k, t}(x)|\lesssim_{N} \frac{2^{2k}}{(1+2^{k}|x-t\gamma|)^N},$ we further derive  
\begin{align}|K^{k}_{j}(x)|\lesssim_{N} 2^{-3j} (2^{k-j}|x|)^{-N},\quad \text{for}\quad |x|\geq c 2^{j}.\label{eq:kernel}\end{align}
Note that for large $|x|$, $|x-t\gamma|\gtrsim |x|$ as $t\in [1,2]$ and $|\gamma|$ is finite.
For $y\in Q$, $z\notin B_x^L$ we have $|y-z|\gtrsim2^L$ and
$|y-z|\ge C_\gamma2^j$ for $j\le L$. Now, by \eqref{eq:kernel}, the global part is bounded by 
\begin{align}
\nonumber&C_{N}\sum_{j\leq L} \frac{1}{|Q|}\int_{Q}\int_{\R^3\setminus B^{L}_{x}} 2^{-3j}\frac{1}{(1+2^{k-j}|y-z|)^N}|f(z)|\, dz\, dy\\
&\lesssim C_{N} 2^{-kN} \|f\|_{L^\infty(\R^3)} \sum_{j\leq L} 2^{-(L-j)(N-3)}\lesssim  C_{N} 2^{-kN} \|f\|_{L^\infty(\R^3)}.\label{glo}
\end{align}
Combining \eqref{loc} and \eqref{glo} we obtain
\begin{align}
\|\mathcal{G}f\|_{L^\infty(\R^3)}\lesssim 2^{-k(\sigma(3)-\epsilon)} \|f\|_{L^\infty(\R^3)}.\label{inf-est}   
\end{align}
 Interpolating \eqref{l3est} and \eqref{inf-est} we get
\begin{equation*}
 \|\mathcal Gf\|_{L^p(\R^3)}\lesssim
 2^{-k(\sigma(3)-\epsilon)}\|f\|_{L^p(\R^3)},\qquad 3<p<\infty.
 \label{eq:G}
\end{equation*}

\medskip

\noindent \textbf{Estimates for $U_n, 1<n<k$:} The triangle inequality in $L^3_t$ followed by the union bound argument gives
\begin{align*}
 U_nf\leq 2\mathfrak{M}\big(\sup_j\|\A^k_jf(\cdot,t)\|_{L^3_t}\big)\lesssim \mathfrak{M}\left(\sum_j\|\A^k_jf\|_{L^3_t}^p\right)^{1/p}.
\end{align*}
Consequently, using the boundedness of $\mathfrak{M},$  union bound argument and local smoothing estimates and finally by Littlewood-Plaey theory we obtain for $p>1$,
\begin{align}
 \|U_nf\|_{L^p(\R^3)}
&\lesssim\left(\sum_j\|\A^k_jf\|_{L^p_x(L^3_t)}^p\right)^{1/p}\lesssim\left(\sum_j\|\A^k_j \widetilde P_{k-j} f\|_{L^p_x(L^p_t)}^p\right)^{1/p}
 \lesssim 2^{-k(\sigma(p)-\epsilon)}\|f\|_{L^p(\R^3)} .
 \label{eq:middle}
\end{align}

\noindent \textbf{Estimates for $U_n, n\geq k$:} Note that, for $Q\in \mathcal{Q}_{L}(x)$ we have
\begin{align}
\frac{1}{|Q|^2}\int_{Q}\int_{Q}\|\widetilde P_{k-L-n}g(y, \cdot)- \widetilde P_{k-L-n}g(w, \cdot)\|_{L^3_{t}}\, dw\, dy \lesssim 2^{k-n} M(\|g\|_{L^3_{t}})(x).
\label{cancel}
\end{align} 
Using this with $g=\A^k_{L+n}f=P_{k-L-n} \A^k_{L+n}f,$ we derive
\begin{align}
\nonumber&\big\|\sup_{L\in\Z}\sup_{Q\in\mathcal Q_L(x)}
 \frac{1}{|Q|^2}\int_Q\int_Q
 \|\A_k^{L+n}f(y,\cdot)-\A_k^{L+n}f(w,\cdot)\|_{L^3_t}\,dw\,dy \big\|_{L^p}\\
 &\nonumber\lesssim 2^{k-n} \|\sup_{j} M(\|\A^k_{j}f\|_{L^3_t})\|_{L^p}\lesssim 2^{k-n} \|\big(\sum_{j} M(\|\A^k_{j}f\|_{L^3_t})^p\big)^{1/p}\|_{L^p}\\
& \lesssim 2^{k-n} \left(\sum_j\|\A^k_j \widetilde P_{k-j} f\|_{L^p_x(L^p_t)}^p\right)^{1/p}
 \lesssim 2^{k-n} 2^{-k(\sigma(p)-\epsilon)}\|f\|_{L^p(\R^3)} \label{large} .
\end{align}

Combining \eqref{eq:split}, \eqref{eq:g}, \eqref{eq:middle} and \eqref{large} we obtain that 
\begin{align*}
\|G\|_{L^p(\R^3)}\lesssim \|G^{\#}f\|_{L^p(\R^3)}&\lesssim \|\mathcal{G}f\|_{L^p(\R^3)}+\sum_{1<n\leq k}\|U_{n}f\|_{L^p(\R^3)}+\sum_{n>k} \|U_{n}f\|_{L^p(\R^3)}\nonumber\\
&\lesssim (1+k) 2^{-k(\sigma(p)-\epsilon)}\|f\|_{L^p(\R^3)}. 
\end{align*}

The preceding
argument therefore applies verbatim to $2^{-k}\partial_t\A^k_j$ and we derive \begin{align*}
 &\left\|\left(\sum_{j\in\Z}\int_1^2
 |A_{2^jt}P_{k-j}f|^3\,dt\right)^{1/3}\right\|_{L^p(\R^3)}\\
 &\quad+2^{-k}\left\|\left(\sum_{j\in\Z}\int_1^2
 |\partial_t[A_{2^jt}P_{k-j}f]|^3\,dt\right)^{1/3}\right\|_{L^p(\R^3)}
 \lesssim  2^{-k(\sigma(p)-\epsilon)}\|f\|_{L^p(\R^3)},
\end{align*}
since $(1+k)2^{-k\epsilon/2}$ is bounded independently of $k$.
\end{proof}

\begin{remark}
In fact, following our arguments, we can prove that for $k\geq 1,$ $3\leq r<p<\infty,$ we have
\begin{align}
&\big\|
\big(
\sum_j\int_1^2
|A_{2^jt}P_{k-j}f|^r\,dt
\big)^{1/r}
\big\|_{L^p(\R^3)}+2^{-k}
\big\|
\big(
\sum_j\int_1^2
\left|\partial_t[A_{2^jt}P_{k-j}f]\right|^r\,dt
\big)^{1/r}\big\|_{L^p(\R^3)}\lesssim c_{p,\epsilon}
2^{-k(\sigma(p)-\epsilon)}\big\|f\big\|_{L^p(\R^3)}.
\label{glo-r}
\end{align}
It will be required for our subsequent purpose.
\end{remark}


\begin{proof}[Proof of Theorem~\ref{no-decay}]
    Note that $\sigma(p)=1/3$ for $3<p\leq 4$. Therefore, the estimate \eqref{ptwise:dyadic:sv}, and Proposition~\ref{global:prop} yields the desired inequality.
\end{proof}

Our next proposition estimates $S_{p}(\mathcal A^k f).$ Recall
$$S_{p}(\mathcal A^k f):=\big(
\sum_j
V_3\big(A_{2^jt}P_{k-j}f(x):t\in[1,2]\big)^p
\big)^{1/p}.$$

\begin{prop}
Let $2 < p < \infty$ and $k \ge 1$. Then
\begin{align}
\|S_p(\mathcal{A}^k f)\|_{L^p(\mathbb{R}^3)} \le C_{p, \epsilon} \, 2^{-k\left(\sigma(p) - \frac{1}{p} - \epsilon\right)} \|f\|_{L^p(\mathbb{R}^3)}.
\label{decay1}
\end{align}    
\label{decay}
\end{prop}

\begin{proof}
Using Sobolev embedding \eqref{sobolev} and Young's inequality we get 
\begin{align*}
S_p(\mathcal A^kf)(x)
\lesssim&
2^{k/p}
\big(
\sum_j\int_1^2
|A_{2^jt}P_{k-j}f(x)|^p\,dt
\big)^{1/p}+
2^{-k(1-1/p)}
\big(
\sum_j\int_1^2
|\partial_t[A_{2^jt}P_{k-j}f(x)]|^p\,dt
\big)^{1/p}.
\end{align*}

The $L^p$ norm of the first term is controlled as 
\begin{align*}
\big(
\sum_j\int_1^2
|A_{2^jt}P_{k-j}f|^p\,dt
\big)^{1/p}\|_{L^p(\R^3)}^p &=
\sum_j\int_1^2
\|A_{2^jt}P_{k-j}f\|_{L^p(\R^3)}^p\,dt\\
&\lesssim c_{p,\epsilon}
2^{-kp(\sigma(p)-\epsilon)}
\sum_j\|\widetilde P_{k-j}f\|_{L^p(\R^3)}^p\lesssim c_{p,\epsilon}
2^{-kp(\sigma(p)-\epsilon)} \|f\|_{L^p(\R^3)}^p.
\end{align*}

Similarly controlling the second term, we obtain
\begin{align*}
\|S_p(\mathcal A^kf)\|_{L^p(\R^3)}
&\lesssim c_{p,\epsilon}
\left[
2^{k/p}
+
2^{-k(1-1/p)}2^k
\right]
2^{-k(\sigma(p)-\epsilon)}\|f\|_{L^p(\R^3)}\\
&\lesssim c_{p,\epsilon}
2^{-k(\sigma(p)-1/p-\epsilon)}\|f\|_{L^p(\R^3)}.
\end{align*}

Observe that $\sigma(p)>\frac{1}{p}$ for all $3<p<\infty$. So we obtain a decay in \eqref{decay1}.
\end{proof}

\subsection{Proof of sufficiency in Theorem~\ref{mainresult}}
We now have enough ingredients to prove the Theorem~\ref{mainresult}. It is divided into several cases. 
\medskip

\begin{proof}[Proof of Theorem~\ref{mainresult}]
We shall divide the proof into three different cases.\\
\textbf{Case: when $3< p\leq4$.} Interpolating the frequency-localized estimates in
Propositions~\ref{no-decay} and~\ref{decay}, we obtain, for every
$q>3$, a constant $\delta=\delta(p,q)>0$ such that
\begin{align*}
\|S_q(\mathcal A^k f)\|_{L^p}
\lesssim  2^{-\delta k}\|f\|_{L^p},
\qquad k\geq1.
\end{align*}
Summing over $k\geq1$ and combining this with the low-frequency
estimate yields
\[
\|S_q(\mathcal A f)\|_{L^p}
\leq C_{p,q}\|f\|_{L^p},
\qquad 3<p\leq4,\quad q>3.
\]
\medskip
\textbf{Case: When $8\leq p<\infty$. }
 Let $q>p/2,$ we need to show that 
    \[\Vert S_{q}(Af)\Vert_{L^{p}}\lesssim \Vert f\Vert_{L^p}.\]

Let us choose $r$ such that $\frac{p}{2}<r<\min\{p,q\}.$ In particular,
$r>4$ and $\frac1r<\frac2p=\sigma(p).$ Fix $\rho$ satisfying
\begin{align*}
\label{rho}
    \frac1r<\rho<\frac2p.
\end{align*}
In this range, for each $k\geq 1,$ the local smoothing estimate gives
\begin{equation*}\label{eq:local-smoothing-large-p}
\left(
    \int_1^2\|A_tP_kf\|_{L^p}^p\,dt
\right)^{1/p}
\lesssim 
2^{-k(2/p-\epsilon)}
\|\widetilde P_kf\|_{L^p}.
\end{equation*}
for every $\epsilon>0.$ Since $r<p$, above implies
$$\big\|\|A_tP_kf\|_{L^r([1,2])}\big\|_{L^p(\R^3)}
    \lesssim
    2^{-k(2/p-\epsilon)}\|f\|_{L^p(\R^3)}.$$
Moreover, the same is true for $2^{-k}\partial_t[A_tP_kf].$ Now using \eqref{glo-r} we derive
\begin{align}
\label{glo-case-2}
&\big\|
\big(
\sum_j\int_1^2
|A_{2^jt}P_{k-j}f|^r\,dt
\big)^{1/r}
\big\|_{L^p(\R^3)}+2^{-k}
\big\|
\big(
\sum_j\int_1^2
\left|\partial_t[A_{2^jt}P_{k-j}f]\right|^r\,dt
\big)^{1/r}\big\|_{L^p(\R^3)}\lesssim c_{p,\epsilon}
2^{-k(\frac{2}{p}-\epsilon)}\big\|f\big\|_{L^p(\R^3)}.
\end{align}

Since we have the pointwise estimate \begin{align*}
S_r(\A^kf)(x)\lesssim
2^{k/r}
\big(
    \sum_j\int_1^2
    |A_{2^jt}P_{k-j}f(x)|^r\,dt
\big)^{1/r}+2^{k/r-k}
\big(
    \sum_j\int_1^2
    |\partial_t[A_{2^jt}P_{k-j}f(x)]|^r\,dt
\big)^{1/r},
\end{align*}
using \eqref{glo-case-2} we obtain
\begin{align*}
    \|S_r(\A^kf)\|_{L^p(\R^3)}
    &\lesssim
    2^{k/r}2^{-k\rho}\|f\|_{L^p(\R^3)}=2^{-k(\rho-1/r)}\|f\|_{L^p(\R^3)},
    \qquad k\geq1.
\end{align*}
Since $\rho>1/r$, summing in $k\geq 1,$ we further derive 
\begin{align*}
\|S_r(\A f)\|_{L^p} \lesssim \|f\|_{L^p}.    
\end{align*}
As $q>r,$ the above implies the desired estimate 
\begin{align*}
\|S_q(\A f)\|_{L^p(\R^3)} \lesssim \|f\|_{L^p(\R^3)}, \quad q>\frac{p}{2}.  
\end{align*}
    
Such a choice is possible because $q>p/2$. \\
\medskip

\noindent \textbf{Case: When $4<p<8$.}
    We want to prove the following estimate 
    \[\Vert S_q(Af)\Vert_{L^p}\lesssim \Vert f\Vert_{L^p},~~~~~\text{for}~~~q>\frac{6p}{p+4}.\]
    The proof of this estimate follows from interpolation. Note that from the above two cases we have got the following estimates 
    \begin{eqnarray*}
        &&\Vert S_{q_0}(Af)\Vert_{L^4}\lesssim \Vert f\Vert_{L^4},~~~\text{for}~~q_0>3,\\
        &&\Vert S_{q_1}(Af)\Vert_{L^8}\lesssim \Vert f\Vert_{L^8},~~~\text{for}~~q_1>4.
    \end{eqnarray*}
    Now, applying interpolation we get that for any $p\in (4,8)$ given by $\frac{1}{p}=\frac{\theta}{4}+\frac{1-\theta}{8}$, we have 
    \[\Vert S_{q}(Af)\Vert_{L^p}\lesssim \Vert f\Vert_{L^p},~~~\text{for}~~\frac{1}{q}=\frac{\theta}{3+\epsilon}+\frac{1-\theta}{4+\epsilon},\]
    where $\epsilon$ is arbitrarily small positive number. Since $\theta=8/p-1$, therefore we get $1/q<(4+p)/6p$.
\end{proof}

\section*{Acknowledgments} AG gratefully acknowledges the 
support by the Industrial Consultancy and Sponsored Research (IC\& SR), Indian Institute of Technology Madras for the
New Faculty Initiation Grant RF25261459MANFIG009296. KS gratefully acknowledges  the vibrant research environment of Kerala School of Mathematics and the support of Inspire Faculty Award of Department of Science
and Technology, Govt. of India, Registration no. IFA23-MA191.

\end{document}